\documentclass[12pt,fleqn]{article}
\usepackage{graphicx}

\usepackage{latexsym}

\usepackage{amsmath}
\usepackage{amsthm}
\usepackage{amssymb}
\usepackage{amsfonts}

\begin{document}

\newcommand{\rf}[1]{(\ref{#1})}
\newcommand{\rff}[2]{(\ref{#1}\ref{#2})}

\newcommand{\ba}{\begin{array}}
\newcommand{\ea}{\end{array}}

\newcommand{\be}{\begin{equation}}
\newcommand{\ee}{\end{equation}}

\newcommand{\const}{{\rm const}}
\newcommand{\ep}{\varepsilon}
\newcommand{\Cl}{{\cal C}}
\newcommand{\rr}{{\vec r}}
\newcommand{\ph}{\varphi}
\newcommand{\R}{{\mathbb R}}  
\newcommand{\C}{{\mathbb C}}  
\newcommand{\T}{{\mathbb T}}  
\newcommand{\Z}{{\mathbb Z}}  
\newcommand{\N}{{\mathbb N}}  

\newcommand{\Eq}{{\cal E}_q}  
\newcommand{\Sq}{{\cal S}in_q} 
\newcommand{\Cq}{{\cal C}os_q}

\newcommand{\e}{{\bf e}}

\newcommand{\m}{\left( \ba{c}}
\newcommand{\ema}{\ea \right)}
\newcommand{\mm}{\left( \ba{cc}}
\newcommand{\miv}{\left( \ba{cccc}}

\newcommand{\av}[1]{\mbox{$\langle #1 \rangle$}} 
\newcommand{\scal}[2]{\mbox{$\langle #1 \! \mid #2 \rangle $}}
\newcommand{\ods}{\par \vspace{0.5cm} \par}
\newcommand{\mods}{\par \vspace{-0.5cm} \par}
\newcommand{\dis}{\displaystyle }
\newcommand{\mc}{\multicolumn}
\newcommand{\id}{{\rm id}}
\newcommand{\no}{\noindent}

\newtheorem{prop}{Proposition}
\newtheorem{Th}[prop]{Theorem} 
\newtheorem{lem}[prop]{Lemma}
\newtheorem{rem}[prop]{Remark}
\newtheorem{cor}[prop]{Corollary}
\newtheorem{Def}[prop]{Definition}
\newtheorem{open}{Open problem}
\newtheorem{ex}{Example}
\newtheorem{exer}{Exercise}

\newenvironment{Proof}{\par \vspace{2ex} \par
\noindent \small {\it Proof:}}{\hfill $\Box$ 
\vspace{2ex} \par }

\title{\bf 
 Improved $q$-exponential and $q$-trigonometric functions}
\author{
 {\bf Jan L.\ Cie\'sli\'nski}\thanks{\footnotesize
 e-mail: \tt janek\,@\,alpha.uwb.edu.pl}
\\ {\footnotesize Uniwersytet w Bia{\l}ymstoku,
Wydzia{\l} Fizyki}
\\ {\footnotesize ul.\ Lipowa 41, 15-424
Bia{\l}ystok, Poland}
}

\date{}

\maketitle

\begin{abstract}
We propose a new definition of the $q$-exponential function. Our $q$-exponential function maps the imaginary axis into the unit circle and the resulting $q$-trigonometric functions are bounded and satisfy the Pythagorean identity. 
\end{abstract}

\ods



\pagebreak

\section{Introduction}

The quantum calculus ($q$-calculus) is an old, classical branch of mathematics, which can be traced back to Euler and Gauss \cite{Er-hist,Ku} with important contributions of Jackson a century ago \cite{Jac1,Jac2}. In recent years there are many new developments and applications of the $q$-calculus in mathematical physics, especially  concerning special functions  \cite{AAR,Ch,ERS,GR} and quantum mechanics \cite{AKW,Bi,Fa,Fi,DP,LRW,Ub}. Many papers were devoted  to various approaches to $q$-deformations of elementary functions, including exponential and trigonometric functions \cite{At-1,At-2,BC,Go,AKK,An,Qu,Su-add,Su-anad}. 

In this paper we propose new definitions of the $q$-exponential function and $q$-trigonometric functions. These results are motivated by recent developments in the time scales calculus, where new exponential, hyperbolic and trigonometric function have been defined \cite{Ci-trig}. 
The concept of time scales unifies difference and differential calculus \cite{Hi}. The $q$-calculus can be considered as a calculus on a special time scale (see, e.g., \cite{BG-Lap}). 

The functions presented in this paper have better qualitative properties than standard $q$-exponential and $q$-trigonometric functions. In order to discuss and compare these properties we begin with a short summary of the classical results, usually following the textbook \cite{KC}. 

In the standard approach to the $q$-calculus two exponential function are used: 
\be  \label{eser}
  e_q^z = \sum_{n=0}^\infty \frac{z^n}{[n]!} \ ,
\qquad E_q^z = \sum_{n=1}^\infty \frac{z^n}{ [\tilde{n}] !}    \ , 
\ee
where $q$ is positive, $z$ is complex, and 
\be  \ba{l} \dis
[n]! = [1][2]\ldots [n] \ , \qquad [k] = 1 + q + q^2 + \ldots + q^{k-1} \ , \\[2ex]\dis
 [\tilde{n}]! = [\tilde 1] [\tilde 2] \ldots [\tilde n] \ , \qquad 
[\tilde{k}] = 1 + \frac{1}{q} + \frac{1}{q^2} + \ldots + \frac{1}{q^{k-1}} \ .
\ea \ee
Hence we immediately get \ $E_q^z = e_{1/q}^z$.   
Another, more popular, form of $E_q^z$ is obtained using the identity  
\be
 [\tilde{n}]! = q^{\frac{(1-n) n}{2}} [n]! \ .
\ee
Both exponential functions can be represented by infinite products,  
\be
e_q^z = \prod_{k=0}^\infty ( 1 - (1-q)q^k z )^{-1} ,  \qquad  
E_q^z =  \prod_{k=0}^\infty ( 1 + (1-q)q^k z )  \ .
\ee 
From this form we easily see that $e_q^z E_q^{-z} = 1$. Moreover, 
\be
 D_q e_q^z = e_q^z \ , \qquad D_q E_q^z = E_q^{q z} \ , 
\ee
where $D_q$ ($q$-derivative or Jackson's derivative) is defined by 
\be
 D_q f (z) := \frac{f (qz) - f (q)}{q z - z}  \ .
\ee

The existence of two representations of $q$-exponential functions (infinite series and infinite product) is related to well known formulae for the usual exponential function ($q=1$), 
\be
  e^z = \sum_{k=0}^\infty \frac{ z^k}{k!} = \lim_{m\rightarrow \infty} \left( 1 + \frac{z}{m} \right)^m \ . 
\ee

Two exponential functions of the quantum calculus generate two 
pairs of the $q$-trigonometric functions. Using notation of \cite{KC} we have:
\be \ba{l}  \dis  \label{scSC} 
\sin_q x = \frac{e_q^{i x} - e_q^{-ix}}{2 i} \ , \quad {\rm Sin}_q x = \frac{E_q^{i x} - E_q^{-ix}}{2 i} \ , \\[3ex] \dis
\cos_q x = \frac{e_q^{i x} + e_q^{-ix}}{2 } \ , \quad {\rm Cos}_q x = \frac{E_q^{i x} + E_q^{-ix}}{2 } \ .
\ea \ee
Taking into account properties of $q$-exponential functions (see above) we easily derive 
properties of standard $q$-trigonometric functions:  
\be  \ba{l} \dis  \label{trig-ident}
 \cos_q x \ {\rm Cos}_q x + \sin_q x \ {\rm Sin}_q x = 1  \  , \\[2ex] \dis
 \sin_q x \ {\rm Cos}_q x = \cos_q x \ {\rm Sin}_q x \ , 
\ea \ee 
\be \ba{l} \dis
  D_q \sin_q x = \cos_q x \ , \quad D_q \cos_q x = - \sin_q x \ , \\[2ex]\dis
D_q {\rm Sin}_q x = {\rm Cos}_q (qx) \ , \quad D_q {\rm Cos}_q x = - {\rm Sin}_q (q x) \ . \\[1ex]
\ea \ee
Note that the corresponding tangents coincide: ${\rm Tan}_q x = \tan_q x$.

\section{Improved $q$-exponential function}

 New $q$-exponential function  ${\cal E}_q^z$ is defined as 
\be  \label{newexp} 
{\cal E}_q^z := e_q^{\frac{z}{2} } \, E_q^{\frac{z}{2}}  = \prod_{k=0}^\infty \frac{1 + q^k (1-q) \frac{z}{2} }{1 - q^k (1-q) \frac{z}{2}} \ ,
\ee 
where $e_q^z$, $E_q^z$ are standard $q$-exponential functions.  
This definition is motivated by the classical Cayley transformation 
\be
   z \rightarrow {\rm cay} (z, a) := \frac{1 + a z}{1 - a z} \ ,
\ee
see, e.g., \cite{Ci-trig,Is-Cay}. Indeed,
\be
      {\cal E}_q^{q z} = \frac{1 - (1-q) \frac{z}{2}}{1 + (1-q) \frac{z}{2} } \  {\cal E}_q^z = {\rm cay} \left( - \frac{z}{2}, \ 1-q \right) \  {\cal E}_q^z \ . 
\ee

\ods

\begin{Th} 
The $q$-exponential function ${\cal E}_q^z$ is analytic in the disc  $|z| < R_q$ and 
\be  \label{Eser}
  {\cal E}_q^z = \sum_{n=0}^\infty  \frac{z^n}{\{n\}!} \ , 
\ee
for $|z| < R_q$, where 
\be \label{Rq}
   R_q = \left\{ \ba{ccl} \dis \frac{2}{1-q} & {\it for} &  0 <   q  < 1 \ , \\[3ex] \dis \frac{2 q}{q -1} & {\it for}  &   q   > 1 \ , 
\\[3ex] \dis \infty & {\it for}  & q = 1 \ , \ea \right. 
\ee
\be \label{paren}
 \{n\} :=  \frac{1+q + \ldots + q^{n-1}}{ \frac{1}{2} (1 + q^{n-1}) } = \frac{[n]}{ \frac{1}{2} (1 + q^{n-1}) }  =  \frac{2 (1-q^n)}{(1-q)(1+q^{n-1})} \ , 
\ee 
\ods
\no and, finally,  $\{n\}! = \{ 1 \} \{ 2 \} \ldots \{ n \}$.
\end{Th}

\begin{Proof} 
In the disc $|z|<1$ both series \rf{eser} are absolutely convergent for any $q \in \R_+$. Multiplying them we get
\be  \label{serprod} 
 e_q^{\frac{z}{2} } \, E_q^{\frac{z}{2}}  = \sum_{k=0}^\infty \sum_{j=0}^\infty  \frac{ q^{ \frac{ j (j-1)}{2} } \left( \frac{z}{2} \right)^{k+j} }{[k]! [j]!} = \sum_{n=0}^\infty \frac{ \left( \frac{z}{2} \right)^n  }{[n]! } \left( \sum_{j=0}^n  \frac{q^{ \frac{ j (j-1)}{2} }  [n]!}{[j]! [n-j]!} \right) \ .
\ee
Using Gauss's binomial formula (see, e.g., \cite{KC}, formula (5.5)) 
\be  
 \prod_{j=0}^{n-1} (z+ a q^{j}) = \sum_{j=0}^n  \frac{q^{ \frac{ j (j-1)}{2} } \ [n]! }{[j]! [n-j]! } \  a^j z^{n-j} \ ,
\ee
we have, as a particular case, 
\be \label{binom} 
 \sum_{j=0}^n \frac{q^{\frac{j (j-1)}{2}} \ [n]! }{ [j]! [n-j]!} = (1+1)(1+q)\ldots (1+q^{n-1}) \ .
\ee
Substituting \rf{binom} into \rf{serprod} we get the formula  \rf{Eser} with $\{ n \}$ defined by \rf{paren}. 
In order to obtain the radius of convergence, we compute
\be
\lim_{n \rightarrow \infty} \left| \frac{z^{n+1}}{\{ n + 1 \}!} \right| \left|   \frac{ \{ n \}! }{z^n} \right| = \lim_{n \rightarrow \infty} \left| \frac{z}{\{ n + 1 \} } \right| =  \left\{ \ba{lll} \frac{(1-q) |z|}{2} & {\rm for} & q < 1 \\[2ex] \frac{(q-1) |z|}{2 q} & {\rm for} & q > 1 \ea \right. 
\ee
Then, using d'Alembert's test, we get (for $q \neq 1$) the radius of convergence \rf{Rq}. Note that $R_{1/q} = R_q$. In the case $q=1$ all $q$-exponential functions coincide with $e^z$, hence $R_1 = \infty$. 
\end{Proof}

\begin{Th} 
The $q$-exponential function ${\cal E}_q^z$ has the following properties:
\be  \label{cechy1} 
  {\cal E}_q^{-z} = \left( {\cal E}_q^z \right)^{-1} \ , \qquad 
| \ {\cal E}_q^{ix} | = 1 \ , 
\ee
\be \label{cechy2}
{\cal E}_q^z = {\cal E}_{1/q}^z \ , \qquad 
 D_q {\cal E}_q^z = \av{ {\cal E}_q^z  } \ , 
\ee
where $z \in \C$, $x \in \R$ and we use the notation $\dis \av{ f(z)} := \frac{f (z) + f (qz) }{2}$. 
\end{Th}

\begin{Proof} The first equation of \rf{cechy1} is a straightforward consequence of the definition \rf{newexp}. Then, 
$\overline{{\cal E}_q^z} = {\cal E}_q^{\bar z}$. Hence, 
\be
   |  {\cal E}_q^{i x} |^2 = \overline{ {\cal E}_q^{i x} }  {\cal E}_q^{i x} =  {\cal E}_q^{- i x}  {\cal E}_q^{i x} = 1 \ . 
\ee
 The symbol $\{ n \}$ depends on $q$. In this proof it is convenient to use more precise notation  $\{ n \} \equiv  \{ n \}_q$, $\{ n \}! \equiv  \{ n \}_q !$.    The  equation ${\cal E}_q^z = {\cal E}_{1/q}^z$ follows immediately from the obvious identity
\be
  \{ n \}_q !  = \{ n \}_{1/q} ! \ . 
\ee
Finally, 
\be
 D_q  {\cal E}_q^z = \frac{  {\cal E}_q^{q z} -  {\cal E}_q^z }{ qz - z } = \frac{ {\cal E}_q^z }{ (q-1) z } \left (  \frac{1 - (1-q) \frac{z}{2}}{1 + (1-q) \frac{z}{2} } - 1  \right) = \frac{ {\cal E}_q^{q z} }{1 + (1-q) \frac{z}{2} } \ , 
\ee
\be
\av{ {\cal E}_q^z } = \frac{1}{2} \left( {\cal E}_q^{qz} + 
{\cal E}_q^z \right) = \frac{1}{2} \left(  \frac{1 - (1-q) \frac{z}{2}}{1 + (1-q) \frac{z}{2} } + 1  \right) {\cal E}_q^z = \frac{  {\cal E}_q^{q z} }{1 + (1-q) \frac{z}{2} } \ ,
\ee
which implies the second equation of \rf{cechy2}. 
\end{Proof}
 
The properties \rf{cechy1} are identical with analogical properties of the exponential function $e^z$. We point out that neither $e_q^z$ nor $E_q^z$ satisfies \rf{cechy1}. Instead, we have $E_q^{-z} e_q^z = 1$.

\section{Improved $q$-trigonometric functions}

New $q$-sine and $q$-cosine functions are defined in a natural way:
\be  \ba{l} \label{newtrig}  \dis
{\cal S}in_q x = \frac{{\cal E}_q^{i x} - { \cal E}_q^{-ix}}{2 i} \ , \quad {\cal C}os_q x = \frac{{\cal E}_q^{i x} + {\cal E}_q^{-ix}}{2 } \ .
\ea \ee
   
\begin{Th}
$q$-Trigonometric functions defined by \rf{newtrig} satisfy:
\be  \ba{l}  \label{cechytrig}
 {\cal C}os_q^2 x + {\cal S}in_q^2 x = 1 \ , \\[2ex]
    D_q {\cal S}in_q x = \av{ {\cal C}os_q x } \ , \\[2ex]
  D_q {\cal C}os_q x = - \av{ {\cal S}in_q x } \ ,
\ea \ee
\end{Th}
\begin{Proof} Properties \rf{cechytrig} follow directly from  \rf{cechy1}, \rf{cechy2} (note that  ${\cal E}_q^{ix} {\cal E}_q^{-ix} = 1$). 
\end{Proof}

\begin{cor}
$q$-Trigonometric functions ${\cal C}os_q x$, ${\cal S}in_q x$ are real for $x \in \R$. Moreover, for any $x \in \R$, we have  
\be  \label{bound} 
-1 \leqslant {\cal C}os_q x  \leqslant 1 \ , \qquad 
-1 \leqslant {\cal S}in_q x  \leqslant 1 \ .
\ee
\end{cor}

\begin{Th} 
New $q$-trigonometric functions can be expressed by standard $q$-trigonometric functions as follows: 
\be  \ba{l} \dis
  {\cal C}os_q 2x = \cos_q x \ {\rm Cos}_q x - \sin_q x \ {\rm Sin}_q x = \frac{1 - \tan_q^2 x}{1 + \tan_q^2 x} \ ,  \\[3ex] \dis 
{\cal S}in_q 2x = \sin_q x \ {\rm Cos}_q x + \cos_q x \ {\rm Sin}_q x = \frac{2 \tan_q x}{1 + \tan_q^2 x} \ .
\ea \ee
\end{Th}
\begin{Proof} First, we compute 
\be \ba{l} \dis 
\cos_q x \ {\rm Cos}_q x - \sin_q x \ {\rm Sin}_q x = 
 \frac{ e_q^{ix} E_q^{ix} + e_q^{-ix} E_q^{-ix} }{2} = {\cal C}os_q 2x \ , \\[3ex] \dis
\sin_q x \ {\rm Cos}_q x + \cos_q x \ {\rm Sin}_q x = 
\frac{ e_q^{ix} E_q^{ix} - e_q^{-ix} E_q^{-ix} }{2} = {\cal S}in_q 2x \ .
\ea \ee
Then, using \rf{trig-ident}, we get
\be \ba{l} \dis 
 {\cal C}os_q 2x = 
\frac{\cos_q x \ {\rm Cos}_q x - \sin_q x \ {\rm Sin}_q x}{\cos_q x \ {\rm Cos}_q x + \sin_q x \ {\rm Sin}_q x} = \frac{1 - \tan_q x \ {\rm Tan}_q x}{1 - \tan_q x \ {\rm Tan}_q x} \ , \\[3ex] \dis 
{\cal S}in_q 2x = 
\frac{\sin_q x \ {\rm Cos}_q x + \cos_q x \ {\rm Sin}_q x}{\cos_q x \ {\rm Cos}_q x + \sin_q x \ {\rm Sin}_q x} = \frac{\tan_q x + {\rm Tan}_q x}{1 - \tan_q x \ {\rm Tan}_q x} \ .
\ea \ee
\no
Taking into account ${\rm Tan}_q x = \tan_q x$ we complete the proof. 
\end{Proof}

\section{Conclusions}

Motivated by the classical Cayley transformation and recent results in the time scales calculus (see \cite{Ci-trig}), we introduced a new definition of the $q$-expo\-nential function. 
Main advantages of the new $q$-exponential function consist in better qualitative properties (i.e., its properties are more similar to properties of $e^z$). In particular, it maps the unitary axis into the unit circle, compare \rf{cechy1}, which implies excellent properties of new trigonometric functions, including formulae \rf{cechytrig} and  boundedness \rf{bound}. 

Especially interesting is the Pythagorean identity: 
$ {\cal C}os_q^2 x + {\cal S}in_q^2 x = 1$. According to our best knowledge, other $q$-defomations of trigonometric functions 
 do not satisfy this property. The same concerns even the  paper \cite{Go}, full of surprising identities.  

Our exponential function is closely related to both popular $q$-exponential functions \rf{eser}. Therefore, proofs and calculations concerning ${\cal E}_q^z$ can be usually done with the help of known results. We plan to express in terms of the new exponential function classical results containing $q$-exponential functions, and we hope to obtain some improvements.


\end{document}